\documentclass[12pt]{iopart}

\def\bI{\mathbf{I}}

\def\bS{\mathbf{S}}
\def\bT{\mathbf{T}}

\def\bY{\mathbf{Y}}

\def\bee{\mathbf{e}}

\def\bi{\mathbf{i}}

\def\bu{\mathbf{u}}
\def\bv{\mathbf{v}}

\def\by{\mathbf{y}}

\def\half{\frac{1}{2}}
\def\mhalf{{-\frac{1}{2}}}
\def\RM{{\rm ROM}}

\usepackage{iopams}
\usepackage{url}
\usepackage{float}
\usepackage{graphicx}
\usepackage{subfig}
\usepackage{adjustbox}
\usepackage{xcolor}
\usepackage{algorithm}
\usepackage{algpseudocode}
\usepackage{pifont}
\usepackage{setspace}

\begin{document}

\title{Solving inverse scattering problems via reduced-order model embedding procedures}

\author{J\"{o}rn Zimmerling$^1$, Vladimir Druskin$^2$, Murthy Guddati$^3$, Elena Cherkaev$^4$, and  Rob~Remis$^5$}

\address{$^1$ Scientific Computing Division, Department of Information Technology, Uppsala University, L\"{a}gerhyddv\"{a}gen~1, 751 05, Uppsala, Sweden}

\address{$^2$ Worcester Polytechnic Institute, 100 Institute Road, Stratton Hall, Worcester, MA 01609, USA}

\address{$^3$ North Carolina State University, 2501 Stinson Drive, Mann Hall 419, Raleigh, NC 27695, USA} 

\address{$^4$ Department of Mathematics, College of Science, University of Utah, 155 South 1400 East, JWB 233, Salt Lake City, UT 84112, USA}

\address{$^5$ Signal Processing Systems Group, Delft University of Technology, Mekelweg 4, 2628 CD Delft, The Netherlands}

\ead{R.F.Remis@tudelft.nl}
\vspace{10pt}
\begin{indented}
\item[]June 2023
\end{indented}

\begin{abstract}
We present a reduced-order model (ROM) methodology for inverse scattering problems in which the reduced-order models are data-driven, i.e. they are constructed directly from data gathered by sensors. Moreover, the entries of the ROM contain localised information about the coefficients of the wave equation. 

We solve the inverse problem by embedding the ROM in physical space. Such an approach is also followed in the theory of ``optimal grids,'' where the ROMs are interpreted as two-point finite-difference discretisations of an underlying set of equations of a first-order continuous system on this special grid. 
Here, we extend this line of work to wave equations and introduce a new embedding technique, which we call \emph{Krein embedding}, since it is inspired by Krein's seminal work on vibrations of a string. In this embedding approach, an adaptive grid and a set of medium parameters can be directly extracted from a ROM and we show that several limitations of optimal grid embeddings can be avoided. Furthermore, we show how Krein embedding is connected to classical optimal grid embedding and that convergence results for optimal grids can be extended to this novel embedding approach. Finally, we also briefly discuss Krein embedding for open domains, that is, semi-infinite domains that extend to infinity in one direction. 
\end{abstract}

%
%
\vspace{2pc}

%
\submitto{\IP}
%
%
%

\section{Introduction}
\label{sec:intro}
In this paper we discuss so-called reduced-order model (ROM) embedding procedures to solve inverse scattering problems. In such a procedure, data-driven ROMs are constructed from spectral impedance data collected at one end of a bounded interval of interest. Subsequently, a ROM is interpreted as a two-point finite-difference discretisation of an underlying set of first-order continuous wave equations. We refer to this step as embedding of the ROM in physical space and it is this embedding procedure that allows us to determine the medium parameters on the spatial interval of interest. 

Two embedding procedures are discussed in this paper, the first, called \emph{optimal grid embedding}, is based on the theory of (optimal) truncated spectral measure grids  \cite{Borcea_etal}. We consider the recovery of a velocity profile from the first poles and residues of the impedance function, the so-called truncated spectral measure. An implementation of this procedure in terms of travel time coordinates is presented in \cite{Borcea_etal2}. Here, we take this optimal grid procedure as a starting point and present optimal grid embedding in terms of spatial coordinates instead of travel time coordinates. 

However, a drawback of optimal grid embedding is that it requires training for a known medium, data, and boundary conditions. More specifically, to retrieve the spatially varying medium parameters on a certain bounded interval of interest, we first have to determine an optimal grid for a homogeneous reference medium. Having this trained reference grid available, the position dependent medium parameters can be determined on the interval of interest. Moreover, the length of this interval must be included in training as well, but this information is not always available. In radar imaging, for example, the distance to the surface may not be known exactly and for problems on semi-infinite domains, the step sizes quickly become very large \cite{Ingerman_etal}.

The use of a reference grid is avoided in the second embedding procedure presented in this paper, which is a new procedure that we call \emph{Krein embedding}, since it is inspired by Krein's seminal work on vibrations of a string \cite{Kac&Krein, Krein&Nudelman}. Furthermore, the length of the interval of interest is not required either. We also show how Krein embedding is related to optimal grid embedding and that convergence results presented in \cite{Borcea_etal} for optimal truncated spectral measure grids essentially carry over to the Krein embedding approach. 

Since a trained reference grid and the length of a reconstruction interval are not required for Krein embedding, such an embedding approach may be applied to scattering problems on semi-infinite domains as well, that is, domains that extend to infinity in one direction. Furthermore, scattering problems on semi-infinite domains characterised by continuous spectral measures can equivalently be described by scattering poles that correspond to passive dissipative systems \cite{Sjostrand&Zworski}. Our approach is then to construct ROMs that can be interpreted as a finite-difference discretisation of a dissipative first-order system. We call this \emph{Krein-Nudelman embedding}, since the case of a dissipative boundary condition was discussed by Krein and Nudelman in \cite{Krein&Nudelman}. 
       
However, for problems on semi-infinite domains a uniqueness problem arises, since the method embeds an impedance function provided in pole-residue form, which cannot distinguish between a lossy bounded domain and a lossless open domain. In the latter case, the spectrum is not a point spectrum but is represented as such. Nevertheless, if we apply a Krein-Nudelman embedding approach to such a problem, we find that the medium profile is actually recovered up until the last reflector, where Krein-Nudelman embedding places an absorbing effective medium to match the (complex) point spectrum of the impedance function. A numerical example will be presented that illustrates this phenomenon and the uniqueness problem is discussed further in the Appendix.

The remainder of this paper is organised as follows. In Section~\ref{sec:eROM} we discuss the construction of the ROMs from spectral data and the optimal grid and Krein embedding procedures that may be used to retrieve the medium parameters on a bounded interval. Subsequently, these two procedures are discussed in detail in Sections~\ref{sec:SLp_Optimal} and \ref{sec:SLp_Krein}, while in Section~\ref{sec:num_ex} a number of numerical examples are presented that illustrate the performance of both embedding procedures. Finally, Krein-Nudelman embedding on a semi-infinite domain is briefly discussed in Section~\ref{sec:absorbing} and the conclusions can be found in Section~\ref{sec:concl}.     

%
%
\section{Embedding of reduced-order models}
\label{sec:eROM}
We are interested in reconstructing the wave speed $c$ in the wave equation from boundary measurements. We formulate the problem in the temporal Laplace domain with complex Laplace frequency $s\in \mathbb{C}$ and reflecting boundary conditions. Specifically, on the bounded interval $[0,L]$ of interest, the governing Laplace-domain equation is given by
\begin{equation}
\label{eq:basic} 
\frac{\mathrm{d}^2u(x,s)}{\mathrm{d}x^2}-  s^2 \frac{1}{c^2(x)}u(x,s)=0, \quad \frac{\mathrm{d}u}{\mathrm{d}x}\Big|_{x=0}=-s,\quad \mbox{and} \quad u|_{x=L}=0.
\end{equation}
This problem is equivalent to the equation for a vibrating string studied by Kac and Krein \cite{Kac&Krein} if $c^{-2}(x)$ is replaced by the mass density of a string. The above equation also follows from the telegrapher equations for a transmission line that is short-circuited at the far-end of the line and a unit current is fed into the near-end of the line. In this case, $c(x)$ represents the  wave speed along the transmission line. 

For a given wave speed profile, equation~(\ref{eq:basic}) essentially represents a regular Sturm-Liouville problem. We are also interested in the corresponding singular case in which equation~(\ref{eq:basic}) is considered on the semi-infinite interval $[0,\infty)$ with the boundary condition for $u$ at $x=L$ replaced by the condition that $|u(x)| \rightarrow 0$ as $x\rightarrow \infty$ for $\mbox{Re}(s)>0$. 

In both the regular and singular case, the associated inverse problem is to recover $c(x)$ from measurements of the impedance function $f(s)=u(0,s)$ for $s$ on some curve in the complex plane. The spectral inverse problem that we consider in this paper consists of the reconstruction of $c(x)$ from poles and residues of the impedance function $f(s)$. In particular, we assume to have access to the first $n$ complex-conjugate pairs of poles $\lambda_{j}$  and residues $y_{j}$ of the function $f(s)$. In an application, this spectral information is not readily available from measurements of $f(s)$. However, it can be retrieved from the measured impedance function $f(s)$ using the vectorfit algorithm \cite{Gustavsen&Semlyen}. In the numerical experiment section, we illustrate this procedure. 

Assuming that we have $n$ complex-conjugate pole-residue pairs available, we can construct the reduced-order model  
\begin{equation}
\label{eq:TF}
f^{\rm ROM}(s)=\sum_{j=1}^{n}  \frac{{y}_{j}}{s+\lambda_{j}} + \frac{\bar{y}_{j}}{s+\bar{\lambda}_{j}}, 
\end{equation}
where the overbar denotes complex conjugation. For a regular Sturm-Liouville problem, the residues $y_{j}$ are real and positive and the $\lambda_{j}$ are purely imaginary. However, for the singular Sturm-Liouville case discussed later this may not be the case and we therefore include conjugation in the above spectral expansion of $f^{\rm ROM}(s)$. Finally, we mention that if we introduce the diagonal matrix
\begin{equation}
\label{eq:Lambdamat}
\bLambda = {\rm diag} (\lambda_1,\lambda_2,...,\lambda_n, \bar{\lambda}_1, \bar{\lambda}_2,...,\bar{\lambda}_n)
\end{equation} 
and the residue vector 
\begin{equation}
\by = 
[
\sqrt{y_{1}},\sqrt{y_{2}},...,\sqrt{y_{n}},\sqrt{\bar{y}_1},\sqrt{\bar{y}_2},...,\sqrt{\bar{y}_n}
]^{T}
\end{equation}
the reduced-order model can also be written as 
\begin{equation}
\label{eq:TF_diag}
f^{\rm ROM}(s)= \by^{T} \left( \bLambda + s \bI \right)^{-1} \by,
\end{equation}
where $\bI$ is the $2n \times 2n$ identity matrix.

\subsection{Building a reduced-order model from spectral data}
\label{subsec:buidlingROM}

\begin{figure}
\centering
\includegraphics[width=0.9\linewidth]{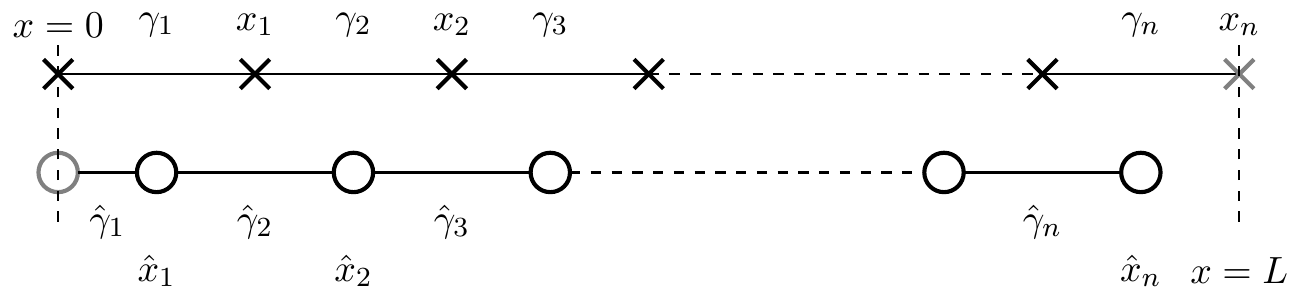}
\caption{Grid used to interpret the ROM as a finite-difference discretisation of the underlying differential operator. The crosses represent primary grid nodes and the circles dual grid nodes. On the grey nodes the boundary conditions are applied.}
\label{fig:grid}
\end{figure}

The key idea behind optimal grid and Krein embedding is to interpret the reduced-order model $f^{\rm ROM}(s)$ as the impedance function of a two-point finite-difference discretisation of an underlying set of first-order differential equations. Specifically, introducing a dual variable $\hat{u}$ and the staggered grid shown in Figure~\ref{fig:grid}, the first-order finite-difference system that corresponds to equation~(\ref{eq:basic}) is given by 
\begin{eqnarray}
\label{eq:GenericFD1}
\frac{\hat u_{j+1}-\hat u_{j}}{\hat\gamma_{j+1}} + s {u}_{j} &= 0, \qquad \forall j=0,\dots,n-1, \nonumber \\
\frac{u_{j}-u_{j-1}}{\gamma_{j}} + s \hat{u}_{j} &= 0, \qquad \forall j=1,\dots,n
\end{eqnarray}
with $u_{j}=u(x_{j})$ and $\hat u_{j}=\hat{u}(\hat{x}_{j})$ and where $\gamma_j$ and $\hat\gamma_j$ are the edge weights, i.e. products of the step sizes $x_{j}-x_{j-1}$ and the medium parameters on the grid. The boundary conditions of the system are given by $\hat{u}(0)=\hat{u}_{0}=1$ and $u(L)=u_{n}=0$. 
 
\noindent
Introducing the $2n$-by-1 vector of unknowns
\begin{equation}
\label{eq:def_vecu}
\bu = 
[u_0,\hat{u}_1,u_1,...,\hat{u}_{n-1},u_{n-1},\hat{u}_n]^T
\end{equation}
the finite-difference system of equation~(\ref{eq:GenericFD1}) can be written compactly as 
\begin{equation}
\label{eq:sysmatvec}
(\bT+s\bI) \bu = \frac{1}{\hat{\gamma}_1} \bee_1, 
\end{equation}
where $\bee_{1}$ is the first canonical basis vector of length $2n$ and matrix $\bT$ is a tridiagonal matrix of order $2n$ given by 
\begin{equation}
\label{eq:Toriginal}
\bT=
\left(
\begin{array}{ccccc}
0 & \hat{\gamma}_1^{-1}  & 0 & &0 \\
-\gamma_1^{-1} & 0 & \gamma_1^{-1} &0& \\
\vdots &  -\hat\gamma_{2}^{-1} & \ddots & \ddots\\
\vdots&     & \ddots & 0 & \hat\gamma_{n}^{-1} \\
0 & \cdots &  0 & -\gamma_{n}^{-1} & 0\\
\end{array}
\right).
\end{equation} 
To accommodate a standard complex-symmetric Lanczos implementation further on, we prefer to work with the   
transpose-symmetric matrix 
\begin{equation}
\label{eq:Tgeneric}
\bT^{\rm TS}=
-\rm{i}
\left(
\begin{array}{ccccc}
0 & (\gamma_1 \hat\gamma_1)^{-\half} & 0 & &0 \\
(\gamma_1 \hat\gamma_1)^{-\half}  & 0 & -(\gamma_{1} \hat\gamma_{2})^{-\half} &0& \\
\vdots &  -(\gamma_{1} \hat\gamma_{2})^{-\half} & \ddots & \ddots\\
\vdots&     & \ddots & 0 & (\gamma_{n} \hat\gamma_{n})^{-\half} \\
0 & \cdots &  0 & (\gamma_{n} \hat\gamma_{n})^{-\half} & 0\\
\end{array}
\right),
\end{equation}
instead of the tridiagonal matrix~$\bT$. Matrix~$\bT^{\rm TS}$ is related to matrix $\bT$ by the similarity transform $\bT^{\rm TS} = \bS \bT \bS^{-1}$, where the similarity matrix $\bS$ is given by
\begin{equation}
\label{eq:symmat}
\bS = {\rm diag}(\sqrt{\hat\gamma_1},\rm{i}\sqrt{\gamma_1},\sqrt{\hat\gamma_{2}},\rm{i}\sqrt{\gamma_{2}},\dots,\sqrt{\gamma_{n}},\rm{i}\sqrt{\gamma_{n}}).
\end{equation}

With $f^{\rm ROM}(s) = u_0(s) = \bee_1^T \bu$ we now find that by solving the system of equations~(\ref{eq:sysmatvec}) for $\bu$ that 
\begin{equation}
\label{eq:TF_tridiag}
f^{\rm ROM}(s)=\frac{1}{\hat\gamma_1} \bee_{1}^{T} (\bT+s\bI)^{-1} \bee_{1} = 
\frac{1}{\hat\gamma_1} \bee_{1}^{T} (\bT^{\rm TS}+s\bI)^{-1} \bee_{1}.
\end{equation}

To extract the finite difference parameters $\gamma_{j},\hat \gamma_{j}$ from the spectral parameters $(\lambda_{j},y_{j})$ we need to transform the diagonal transfer function representation of equation~(\ref{eq:TF_diag}) to the tridiagonal one of equation~(\ref{eq:TF_tridiag}). To that end, let $\bY$ denote the eigenvector matrix of $\bT^{\rm TS}$ such that it satisfies $\bT^{\rm TS} \bY = \bY \bSigma$ with $\bY^{T} \bY = \bI$ and $\bSigma$ is a diagonal matrix of order $2n$ with the eigenvalues of $\bT^{\rm TS}$ on the diagonal. The impedance function can now be written as 	
\begin{equation}	
\label{eq:ROM_TF}	
f^{\rm ROM}(s)=	
\frac{1}{\hat\gamma_{1}} 	
\by_{1}^{T} (\bSigma+s\bI)^{-1} \by_{1},	
\end{equation}	
where $\by_1=\bY^{T}\bee_1$ is the vector containing the first components of all eigenvectors. 	
Setting $\by_1=\sqrt{\hat{\gamma}_1} \by$ and with $\bSigma=\bLambda$, we observe that the above reduced-order model coincides with the reduced-order model of equation~(\ref{eq:TF_diag}). 	
Since the residue vector $\by$ and the poles $\bLambda$ are known, the problem turns into an inverse eigenvalue problem in which we attempt to reconstruct the tridiagonal matrix $\bT^{\rm TS}$ from the first components of its eigenvectors and its known eigenvalues. As is well known, this problem can be solved via the Lanczos algorithm. Specifically, with orthogonalization in the transpose bilinear form \cite{Marshall,Saad}, vector $\by/\sqrt{\by^{T}\by}$ as a starting vector, and the matrix of eigenvalues $\bLambda$ as iteration matrix, the Lanczos algorithm produces the desired tridiagonal matrix. The algorithm is given in Algorithm~\ref{alg:algo1} and as soon as the tridiagonal matrix is obtained, the coefficients $\gamma_{j}$ and $\hat\gamma_{j}$ can be extracted using the recursive scheme shown in Algorithm~\ref{alg:algo2}. We note that the weight $\hat{\gamma}_1$ can be determined directly from the residues, since 	
\begin{equation}
\label{eq:gamma1}	
\by^{T} \by = \hat{\gamma}_1^{-1} \by_1^T \by_1 = \hat{\gamma}_1^{-1}.	
\end{equation}
Finally, we also note that the algorithm explicitly computes the diagonal elements $\alpha_j$ of matrix $\bT^{\rm TS}$, while we know that these elements vanish for lossless media. However, this very same algorithm can be used for problems on open domains or for problems involving lossy media for which the diagonal elements of matrix $\bT^{\rm TS}$ may no longer vanish. Therefore, we prefer to work with this Lanczos algorithm since it accommodates all cases of interest. 

%
%
\smallskip
\noindent
\begin{algorithm}
\caption{The Lanczos algorithm for complex symmetric matrices to obtain the tridiagonal matrix $\bT^{\rm TS}$ from the poles $\bLambda$ and residues $\by$. In this algorithm, $\bY_j$ indicates the $j$th column of the eigenvector matrix $\bY$, $\alpha_{j}$ is the $j$th element on the diagonal $(\bT^{\rm TS})_{j,j}$ and $\beta_{j}$ the $j$th element on the super diagonal $(\bT^{\rm TS})_{j-1,j}$ of matrix $\bT^{\rm TS}$.  }
\label{alg:algo1}
\begin{spacing}{0.9}
\begin{algorithmic}[1]
\Procedure{Complex Symmetric Lanczos}{}
\State Normalize: $\bY_1 = \by/\sqrt{\by^{T} \by}$
\State $\bv_1 = \bLambda\bY_1$ 
\State $\alpha_1 = \bv_1^T \bY_1$
\State $\bv_1 = \bv_1 - \alpha_1 \bY_1$
\For{$j = 2, \ldots, 2n$ }
\State $ \beta_j = \sqrt{\bv_{j-1}^T \bv_{j-1}}$ 
\If{$\beta_j \neq 0$}
\State $\bY_j = \bv_{j-1}/\beta_j$
\State $\bv_j = \bLambda \bY_j$ 
\State $\alpha_j = \bv_j^T \bY_j$ 
\State $ \bv_j = \bv_j - \alpha_j \bY_j - \beta_j \bY_{j-1}$
\Else \,\,Breakdown
\EndIf
\EndFor
\EndProcedure
\end{algorithmic}
\end{spacing}
\end{algorithm}

\begin{algorithm}
\caption{Extract the ROM parameters $\gamma_j$ and $\hat\gamma_j$ from $\beta_{j}$ the super diagonal of the matrix $\bT^{\rm TS}$}
\label{alg:algo2}
\begin{spacing}{0.9}
\begin{algorithmic}[1]
\Procedure{Extract Gamma}{}
\State $\hat{\gamma}_{1}=(\by_1^{T} \by_1)^{{-1}}$
\For{$j = 1, \ldots, n-1$ }
\State $\gamma_j = - ({\beta_{2j}^2}{\hat \gamma_j})^{{-1}}$
\State $\hat \gamma_{j+1} = - (\beta_{2j+1}^2 \gamma_j)^{-1}$
\EndFor
\State $\gamma_n = - ({\beta_{2n}^2}{\hat \gamma_n})^{{-1}}$
\EndProcedure
\end{algorithmic}
\end{spacing}
\end{algorithm}

\subsection{Krein and optimal grid embedding}

Up to this point we constructed a ROM from the poles and residues of a boundary impedance function. The ROM can be interpreted as the impedance function of a two-point finite-difference discretization of a first-order system with primary and dual grid coefficients $\gamma_{j}$ and $\hat\gamma_{j}$, respectively. How to interpret these coefficients depends on the dual variable $\hat{u}$ that is introduced to obtain a first-order system of ODEs from the original second-order ODE. In particular, introducing the dual variable $\hat{u} = -s^{-1} \partial_x u$,  the second-order system of equation~(\ref{eq:basic}) can be written as 

\begin{equation}
\label{eq:KreinCForm}
\left(
\begin{array}{cc}
s&c^2\partial_x \\
\partial_x & s \\
\end{array}
\right)
\left(
\begin{array}{c}
u(x,s) \\
\hat{u}(x,s) 
\end{array}
\right)
=0,
\quad
u(L,s)=0, \quad \hat{u}(0,s)=1,
\end{equation}
and this form will lead to what we call the \emph{Krein embedding} interpretation. On the other hand, introducing the quantities $w=c^{-\half}u$ and $\hat w=c^{\half}\hat u$ as well as the slowness coordinates (sometimes called travel time coordinates)
\begin{equation}
\label{eq:ttc}
T(x)=\int_0^x \frac{1}{c(\xi)}\, \rm{d} \xi
\end{equation}
we obtain the first-order system
\begin{equation}
\label{eq:fo_ttc}
\left(
\begin{array}{cc}
s &c^{\half} \partial_T c^{\mhalf} \\
c^{\mhalf} \partial_T c^{\half}& s
\end{array}
\right)
\left(
\begin{array}{c}
w(T,s) \\
\hat w(T,s)
\end{array}
\right)
=0,
\end{equation}
with $w(T(L),s)=0$ and $\hat{w}(0,s)=c^{\half}(0)$, and this system will leads to a standard optimal grid embedding interpretation. Note that in this latter case the impedance function is $f(s)=u(0,s)=c^{\half}(0) w(0,s)$ and that standard optimal grid embedding is formulated in terms of travel time coordinates. 

In the next two sections we will give detailed descriptions of optimal grid and Krein embedding procedures for wave propagation problems. Convergence of these procedures is also briefly discussed.

\section{Optimal grid embedding}
\label{sec:SLp_Optimal}

Optimal grid embedding of ROMs was developed for the diffusion equation in \cite{Borcea_etal}. In this section, we extend this embedding approach to wave propagation problems and rely on results obtained for the diffusion equation. 

The main difficulty in ROM embedding is that each of the finite difference weights $\gamma_{j}$ and  $\hat\gamma_{j}$ consist of a product of the unknown local medium parameter and an unknown grid step. Fortunately, for diffusion problems it has been shown that there exists a computable grid $\{\gamma^{0}_{j}, \hat\gamma^{0}_{j}\}$ that is independent of the medium parameters in the limit $n\to \infty$ (Lemma 3.2 of \cite{Borcea_etal}). The existence of this grid allows us to obtain the local material parameters by taking ratios of the grid steps $\gamma^{0}_{j}, \hat\gamma^{0}_{j}$ and ROM parameters $\gamma_{j}, \hat\gamma_{j}$ thereby reconstructing the medium. Furthermore, it can be shown the medium estimates converge pointwise in $L^{1}(0,L)$ (Theorem 6.1 of \cite{Borcea_etal}).

As a first step towards optimal grid embedding for wave propagation, we extend the results from \cite{Borcea_etal} by considering equation~(\ref{eq:fo_ttc}) in spatial coordinates instead of travel time coordinates. Subsequently, we take the kinematic effects into account and discuss medium parameter retrieval based on equation~(\ref{eq:fo_ttc}).  

To avoid confusion, we call the medium parameter $\zeta$ when considering (\ref{eq:fo_ttc}) in spatial coordinates. In other words, we start by considering 
\begin{equation}
\label{eq:KreinOp2}
\left(
\begin{array}{cc}
s &\zeta^{\half} \partial_x\zeta^{\mhalf} \\
\zeta^{\mhalf} \partial_x \zeta^{\half}& s
\end{array}
\right)
\left(
\begin{array}{c}
v(x,s)\\
\hat v(x,s)
\end{array}
\right)
=0
\end{equation}
with $v(L,s)=0$, $\hat{v}(0,s)=\zeta^{\half}(0)$ and transfer function $f(s)=\zeta^{\half}(0)v(0,s)$. A staggered finite-difference grid with primary grid steps $\delta_j$ and dual grid steps $\hat \delta_j$ has the grid points
\begin{equation}
x_j = \sum_{i=1}^j \delta_i \quad \mbox{and} \quad \hat x_j = \sum_{i=1}^j \hat\delta_i
\end{equation}
and a two-point finite-difference discretisation of the first-order system (\ref{eq:KreinOp2}) gives
\begin{eqnarray}
\frac{ \hat\zeta_{j+1}^{\mhalf} \hat v_{j+1}-\hat\zeta_j^{\mhalf} \hat v_{j}}{\zeta_{j}^\half\hat\delta_{j+1}} + s {v}_{j} &= 0 \quad \forall j=0,\dots,n-1 \nonumber \\
\label{eq:GenericFD1_2}
\frac{\zeta_j^\half v_{j}- \zeta_{j-1}^\half  v_{j-1}}{\hat\zeta_j^\mhalf \delta_{j}} + s \hat{v}_{j} &= 0 \quad \forall j=1,\dots,n,
\end{eqnarray}
where the notation $v_j=v(x_j,s)$ and $\hat v_j=\hat v(\hat x_j,s)$ is used. Furthermore, $\zeta_j$ and $\hat \zeta_j$ denote averaged medium parameters at the grid locations $x_i$ and $\hat x_i$, respectively. After symmetrization with diagonal matrices the finite-difference pencil that corresponds to (\ref{eq:KreinOp2}) can be written in terms of the transpose-symmetric tridiagonal matrix
\begin{eqnarray}
\label{eq:Toptg}
&\hspace{-2.5cm}\bT^{\rm{TS;og}} =  -\mbox{i}
 \nonumber \\
&\hspace{-2.5cm}\left(
\begin{array}{ccccc}
0 & (\zeta_1\delta_{1}  \hat\zeta_1^{-1}\hat\delta_{1})^{-\half} & 0 & &0 \\
(\zeta_1\delta_{1} \hat\zeta_1^{-1}\hat\delta_{1})^{-\half}  & 0 & -(\zeta_1\delta_{1} \hat\zeta_2^{-1}\hat\delta_{2})^{-\half} &0& \\
\vdots &  -(\zeta_1\delta_{1}  \hat\zeta_2^{-1}\hat\delta_{2})^{-\half} & \ddots & \ddots\\
\vdots&     & \ddots & 0 & (\zeta_n\delta_{n} \hat\zeta_n^{-1}\hat\delta_{n})^{-\half} \\
0 & \cdots &  0 & (\zeta_n\delta_{n} \hat\zeta_n^{-1}\hat\delta_{n})^{-\half} & 0\\
\end{array} 
\right). \nonumber \\
\end{eqnarray}

In ROM embedding we try to interpret the ROM constructed from the measurement data as a discretisation of the underlying equation. If we compare the above tridiagonal discretisation stencil to the tridiagonal ROM matrix from equation~(\ref{eq:Tgeneric}) we find that the discretisation has twice as many unknowns as the ROM has parameters, i.e. we can't disentangle the grid steps from the local medium parameters.

In \cite{Borcea_etal} it was shown that a tridiagonal ROM that matches the lowest $2n$ poles and residues of the transfer function corresponds to a discretisation on a special grid, also known as the optimal grid or spectrally matched grid. This grid is independent of the medium parameter $\zeta$ in the asymptotic limit $n\to \infty$ and can be computed from the ROM of a reference simulation with $\zeta_0(x)=1$. Let this reference grid be characterised by the primary and dual weights $\gamma_{j}^0$ and $\hat \gamma_j^{0}$, respectively, then pointwise estimates of $\zeta$ can be directly extracted from the ROM. To be more specific, let $\zeta^{\RM}(x)$ and $\hat{\zeta}^{\RM}(x)$ be interpolants with interpolation properties
\begin{equation}
\label{eq:impOPTgrid1}
\zeta^\RM(\hat x^{\rm{optimal}}_j)=\gamma_j/\gamma_j^0 \quad \mbox{where} \quad \hat x^{\rm{optimal}}_j=\sum_{i=1}^j \hat\gamma_i^0 
\end{equation}
and
\begin{equation}
\hat\zeta^\RM(x^{\rm{optimal}}_{j-1})=\hat\gamma_j^0/\hat\gamma_j 	\quad \mbox{where} \quad       x^{\rm{optimal}}_j=\sum_{i=1}^j \gamma_i^0, 
\label{eq:impOPTgrid2}
\end{equation}
with $x^{\rm{optimal}}_0=0$, then it can be shown that $\zeta^\RM$ and $\hat\zeta^\RM$ converge pointwise in $L^1(0,L)$ to the true medium profile $\zeta$ as $n\to\infty$, see \cite{Borcea_etal} for details.

Having discussed equation~(\ref{eq:KreinOp2}) using the results of \cite{Borcea_etal}, let us now include kinematic wave effects and consider equation~(\ref{eq:fo_ttc}) to obtain an optimal grid reconstruction scheme for wave propagation problems. Initially following the same procedure as above, we can show that the ratios $\gamma_j/\gamma_j^0 $ and $\hat\gamma_j^0/\hat\gamma_j$ converge to the wave speed $c[x(T)]$ parametrised in slowness coordinates and the optimal grids in equations (\ref{eq:impOPTgrid1}) and (\ref{eq:impOPTgrid2}) are primary and dual optimal slowness grids $T_j^{\rm{opt}}$ and $\hat T_j^{\rm{opt}}$, respectively. To obtain the wave speed in physical coordinates, the inverse slowness transform $x(T):T \mapsto x$ given by
\begin{equation}
\label{eq:xT}
x(T)=\int_0^T c[x(\tau)] \,\mbox{d} \tau
\end{equation}
needs to be extracted as well. This can be realized in two steps: First the optimal grid is adjusted to the average slowness of the medium and, second, the grid is locally adjusted to slowness coordinates. 

The background grid should be computed for the domain $[0,T(L)]$, however this requires knowledge of the average slowness
\begin{equation}
\label{eq:ooc}
\overline{c^{-1}}=\frac{1}{L} \int_0^L c^{-1}(x)\, \mbox{d}x
\end{equation}
as this defines $T(L)=L\overline{c^{-1}}$. There are many ways to extract the average slowness from the ROM and we choose to extracted it from the first ROM coefficient $\hat\gamma_1$ as
\begin{equation}
\label{eq:ooca}
\overline{c^{-1}}\approx \frac{\hat{\gamma}^0_1}{\hat{\gamma}_1 c(0)}, 
\end{equation}
for which obviously the wave speed $c(0)$ at the sensor location is required. The effectiveness of the above equation is due to the fact that $\hat\gamma_1^0$ depends linearly on the domain size and in the limit $n\to \infty$ the ratio $\frac{\hat{\gamma}^0_1}{\hat{\gamma}_1 \overline{c^{-1}}}$ converges to the wave speed at $x=0$. Alternatively, $T(L)$ can be extracted from the limit $\lim_{s\to 0} \frac{1}{s} f(s)$.

With this in place the wave speed $c(x)$ can be extracted from the ROM parameters as
\begin{eqnarray}
\hat{c}^{\rm ROM}(\hat{x}^{\rm{opt}}_j) &= \frac{1}{\overline{c^{-1}}}\,\frac{\gamma_j}{\gamma^0_j} \quad \forall j=1,\dots,n\\
{c}^{\rm ROM}({x}^{\rm{opt}}_{j-1}) & = \overline{c^{-1}} \,\frac{\hat{\gamma}^0_j}{\hat{\gamma}_j} \quad \forall j=1,\dots,n,
\end{eqnarray}
which are to be interpreted as pointwise estimates of the wave speed at the optimal grid points that are adjusted to the local slowness coordinates
\begin{equation}
\label{eq:xopt}
{x}^{\rm{opt}}_j= x(T_j^{\rm{opt}}) = \overline{c^{-1}} 
\sum_{k=1}^j {\gamma}_{k}^0 \hat{c}^{\rm ROM}[x(\hat T^{\rm{opt}}_j)]
\end{equation}
and 
\begin{equation}
\label{eq:xhopt}
\hat{x}^{\rm{opt}}_j=x(\hat T_j^{\rm{opt}})=\overline{c^{-1}} 
\sum_{k=1}^j \hat{\gamma}_{k}^0 {c}^{\rm ROM} [x(T^{\rm{opt}}_j)].
\end{equation}
These estimates converge pointwise in $L^1(0,T(L))$ to the true wave speed which can be obtained as a corollary to Theorem 6.1 in \cite{Borcea_etal}. Essentially, the optimal grid embedding recovers $c[x(T)]$ at points on the optimal grid in slowness coordinates. This recovered $c[x(T)]$ then provides the inverse slowness transform transform $T \mapsto x$, to embed $c(x)$ into physical space.

\section{Krein embedding}
\label{sec:SLp_Krein}

From regular Sturm-Liouville theory it is well known that equation~(\ref{eq:basic})  is satisfied for infinitely many eigenpairs ($\phi_i(x),\lambda_i$), with real eigenfunctions $\phi_i(x)$ and  imaginary eigenvalues $s=\lambda_i$. The true impedance function 
\begin{equation}
f(s)= \sum_{i=1}^{\infty} \frac{y_i}{s+\lambda_i} +\frac{\bar{y}_i}{s+\bar\lambda_i}
\end{equation}
is thus a meromorphic function with an infinite number of poles, corresponding to the eigenvalues $\lambda_i$ and $\bar{\lambda}_i$. Now Krein and optimal grid embedding approaches both utilise truncated spectral impedance data, and in \cite{Kac&Krein} it was shown by Kac and Krein that there is a one-to-one correspondence between an $n$-term truncated spectral impedance function and a medium with a mass function $M(x)$ with $n$ points of increase. Krein embedding allows us to link this mass function  to a continuous wave speed profile $c(x)$.  The following proposition provides an explicit connection between the continuous wave equation and the discrete finite-difference problem obtained from truncated spectral impedance data.

{\sc Proposition~1} [Kac and Krein]
Let  $w_n$ satisfy the wave equation with the mass function $M_n$
\begin{equation}
\label{eq:1Dn} 
\frac{\mathrm{d}^2 w_n(x,s)}{\mathrm{d}x^2} -s^2 M_n(x) w_n(x,s)=0 \quad \frac{\mathrm{d}w_n}{\mathrm{d}x}|_{x=0}=-s,\  w_n|_{x=x_{n+1}}=0,
\end{equation}
in the weak sense, where 
\begin{equation}
\label{eq:mapprox} 
\hspace{-2cm} M_n(x)=\sum_{i=0}^{n-1}\hat\gamma_{i+1} \delta (x-x_i), \quad x_0=0, 
\quad {\rm and} \quad 
x_i=\sum_{k=1}^{i}
\gamma_k \quad {\rm for} \quad  i=1,2,..,n.
\end{equation}
The continuous function $w_n$ then interpolates the finite-difference approximation from equation~(\ref{eq:GenericFD1}), that is, we have $w_n(x_i)=u_i$ for $i=0,1,...,n$.  

\bigskip
\noindent
\textit{Proof:} 
The proof is straightforwardly obtained by substituting a (piecewise) linear interpolation of $u_i$ between the $x_i$'s into the second-order finite-difference equation.
More precisely, let $\tilde{u}(x)$ be the linear interpolation of $u_{i}$ on the grid $\{x_{i}\}$. For $\tilde{u}(x)$ we have 
\begin{eqnarray}
\tilde{u}(x_{i})&=u_{i}, \\
 \frac{{\rm d}}{{\rm d}x} \tilde{u}|_{x=(x_{i}+\delta)}&=\frac{u_{i+1}-u_{i}}{\gamma_{i+1}}, \qquad  0<\delta<{\gamma}_{i+1}\\
\frac{{\rm d}}{{\rm d}x}\tilde{u}|_{x=(x_{i}-\delta)}&=\frac{u_{i}-u_{i-1}}{\gamma_{i}}, \qquad 0<\delta<{\gamma}_{i}.
\end{eqnarray}
Substituting this into the second-order form of the finite-difference relation from equation~(\ref{eq:GenericFD1}) yields
\begin{equation}
\hspace{-1cm} \frac{{\rm d}}{{\rm d}x}\tilde{u}|_{x=(x_{i}+\delta)} -\frac{{\rm d}}{{\rm d}x}\tilde{u}|_{x=(x_{i}-\delta)} - s^{2} u(x_{i}) \hat{\gamma}_{i+1} =0\quad  { \rm for } \quad 0<\delta<\mbox{min}({\gamma}_{i},{\gamma}_{i+1})
\end{equation}
and after integration, we obtain the weak form of equation (\ref{eq:1Dn}) 
\begin{equation}
\int_{x_i-\delta}^{x_i+\delta}
\frac{\mathrm{d}^2}{\mathrm{d}x^2}\tilde{u}(x) {\rm d}x - s^{2} \int_{x_i-\delta}^{x_i+\delta}\tilde{u}(x) \delta(x-x_{i})\hat{\gamma}_{i+1} {\rm d}x =0. \qquad\qquad \square
\end{equation}

\bigskip
\noindent
In Krein embedding we interpret the ROM as the impedance function at $x=0$ of a finite-difference discretization of the first order system~(\ref{eq:KreinCForm}). In particular, a two-point finite-difference discretization on a staggered grid with primary stepsizes $\delta_{j}$ and dual stepsizes $\hat\delta_{j}$ reads
\begin{eqnarray}
\label{eq:GenericFD1_3}
\frac{\hat u_{j+1}- \hat u_{j}}{c_{j+1}^{-2} \hat\delta_{j+1}} + s {u}_{j} &= 0 \quad \forall j=0,\dots,n-1 \nonumber\\
\frac{u_{j}-  u_{j-1}}{\delta_{j}} + s \hat{u}_{j} &= 0 \quad \forall j=1,\dots,n,
\end{eqnarray}
where we used the shorthand notation $u_{j}=u(x_{j},s)$ again. The primary ROM edge-weights $\gamma_{j}$ are therefore interpreted as a step size $\delta_j$ and the dual edge-weights $\hat\gamma_{j}$ as $c_{j}^{-2}\hat\delta_{j}$. To be more precise, we introduce the nondecreasing mass function 
\begin{equation}
\label{eq:mass}
M(x)=\int_{\xi=0}^{x} c^{-2}(\xi)\, {\rm d} \xi
\end{equation}
and the coordinate transform $T \mapsto x$ given by 
\begin{equation}
\label{eq:ct}
x(T)=\int_{\tau=0}^{T} c[x(\tau)] \, {\rm d} \tau.
\end{equation}
The mass function can now be written as  
\begin{equation}
\label{eq:Mf_sc}
M[x(T)]=\int_{\tau=0}^{T} \frac{1}{c[x(\tau)]}\, {\rm d} \tau
\end{equation}
and we observe that 
\begin{equation}
\label{eq:approx}
x^{\rm Krein}_{0}=0, \quad 
x^{\rm Krein}_{j}=\sum_{k=1}^{j} \gamma_{k}, 
{\quad \mbox{and} \quad 
 M_{n}(x^{\rm Krein}_{{ j-1}}) = \sum_{k=1}^{j}\hat\gamma_{k},}
\end{equation}
{ for $j=1,2,...,n$} can be interpreted as ROM quadrature rules of (\ref{eq:ct}) and (\ref{eq:Mf_sc}), respectively. { Finally, if we let $T_j^{\rm Krein}$ denote the slowness coordinate corresponding to $x^{\rm Krein}_{j}$} then, using the results of the optimal grid case discussed in the previous section, we have  
\begin{equation}
\label{eq:SL_KRein_TF}
\chi_n[T_j^{\rm Krein}]= x_{j}^{\rm Krein}=
\sum_{k=1}^{j} \gamma_{k}=
\overline{c^{-1}}
\sum_{k=1}^{j} c^{\rm ROM}(\hat{x}_{k}^{\rm opt})  \, \gamma_{k}^{0}
\end{equation}
and 
\begin{equation}
\label{eq:SL_KRein_M}
M_{n}[\chi_n(T^{\rm Krein}_{ j-1})] = \sum_{k=1}^{j} \hat{\gamma}_{ k} = 
 \overline{c^{-1}}
 \sum_{k=1}^{j}  c^{\rm ROM}(x_{k-1}^{opt})^{-1} \hat{\gamma}^{0}_k.
\end{equation}
Convergence of the Krein embedding follows directly from equations~(\ref{eq:SL_KRein_TF}) and (\ref{eq:SL_KRein_M}). Since the optimal grid parameters converge, convergence of $\chi_n(T)\mapsto x(T)$ and $M_n(x)\mapsto M(x)$ follows again as a corollary of Theorem 6.1 in \cite{Borcea_etal}.

\section{Numerical examples}
\label{sec:num_ex}

Optimal grid embedding leads to pointwise estimates for the wave speed $c(x)$, whereas Krein embedding leads to an estimation of the mass function $M(x)=\int_0^x c^{-2}(\xi) {\rm d} \xi$. Therefore, it is natural to display the inversion result of the optimal grid embedding and Krein embedding in terms of these two quantities. 

\begin{figure}
\centering
\includegraphics[width=0.7\linewidth, trim={0 5cm 0 5cm},clip]{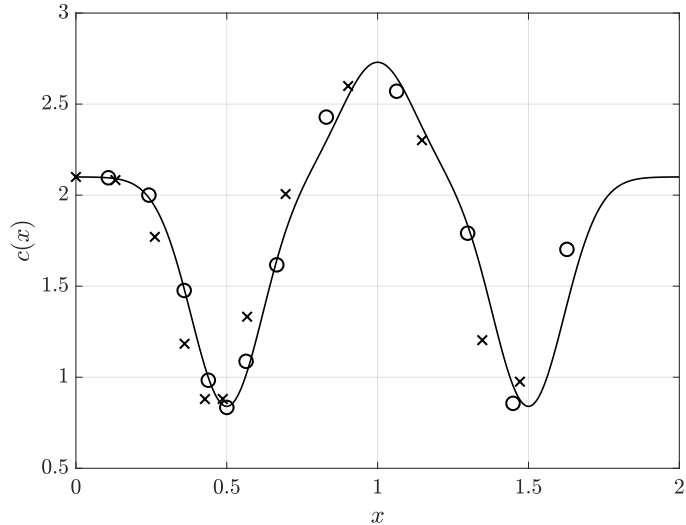}
\caption{Velocity profile $c(x)$ (solid line) and optimal grid reconstruction of this profile based on $n=12$ pole-residue pairs of the impedance function at $x=0$. Crosses: reconstructed velocity values at primary nodes. Circles: reconstructed velocity values at dual nodes.} 
\label{fig:vel_profile1_og}
\end{figure}

In Figure~\ref{fig:vel_profile1_og} the dashed line signifies a smoothly varying velocity profile on the interval $[0,2]$ that we attempt to reconstruct using knowledge of poles and residues of the impedance function at $x=0$. This spectral data is obtained by applying the vectorfit algorithm~\cite{Gustavsen&Semlyen} to the impedance function at $x=0$ to obtain a rational function representation in pole-residue form as indicated in equation~(\ref{eq:TF}).    

\begin{figure}
\centering 
\includegraphics[width=0.7\linewidth, trim={2cm 7cm 2cm 7cm},clip]{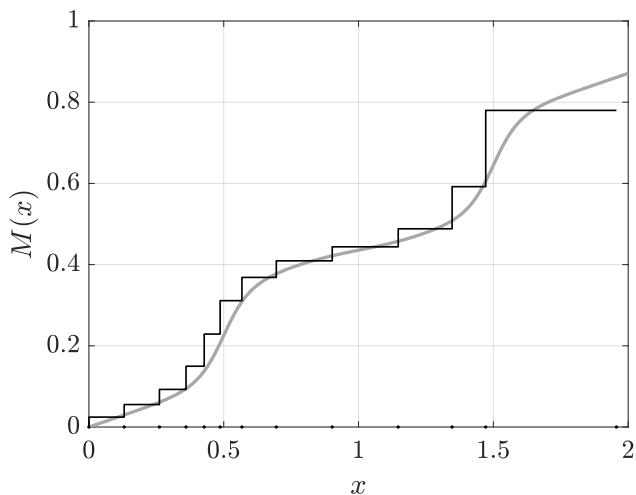}
\caption{Mass function $M(x)$ (solid grey line) and Krein embedding reconstruction of this function (solid black line) based on $n=12$ pole-residue pairs of the impedance function at $x=0$. The Krein grid is visualized with small dots on the x-axis. } 
\label{fig:vel_profile1_K}
\end{figure}

We start with $n=12$ pole-residue pairs of the impedance function $f(s)$ to construct a ROM and realise optimal grid and Krein embeddings. In this example, the background grid in case of optimal grid embedding is extracted from a reference simulation with a constant $c_0=1$. With only $n=12$ spectral points, the wave speed $c(x)$ in case of an optimal grid embedding and the mass function $M(x)$ in case of Krein embedding are accurately reconstructed as illustrated in Figures~\ref{fig:vel_profile1_og} and \ref{fig:vel_profile1_K}, respectively. We note that due to the low number of pole-residue pairs that are used, the estimates are slightly misplaced as the map $T\mapsto x$ is estimated with only a few quadrature points. 


\begin{figure}
\centering
\includegraphics[width=0.7\linewidth, trim={0 5cm 0 5cm},clip]{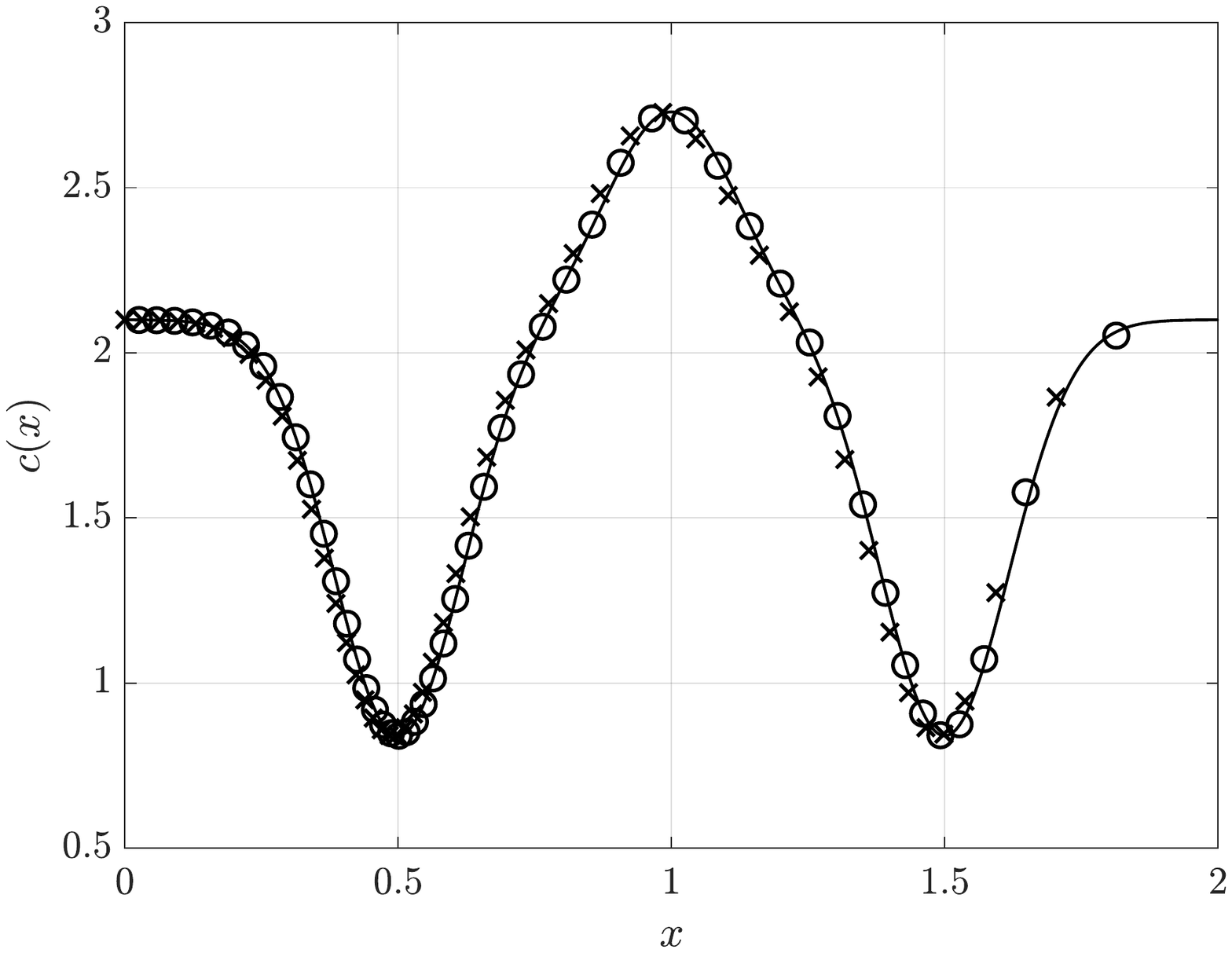}
\caption{Optimal grid embedding reconstruction of the velocity profile (solid line) based on $n=50$ pole-residue pairs of the impedance function at $x=0$. Crosses: reconstructed velocity values at primary nodes, circles: reconstructed velocity values at dual nodes.}
\label{fig:vel_profile1_og_mp}
\end{figure}
\begin{figure}
\centering
\includegraphics[width=0.7\linewidth, trim={0 5cm 0 5cm},clip]{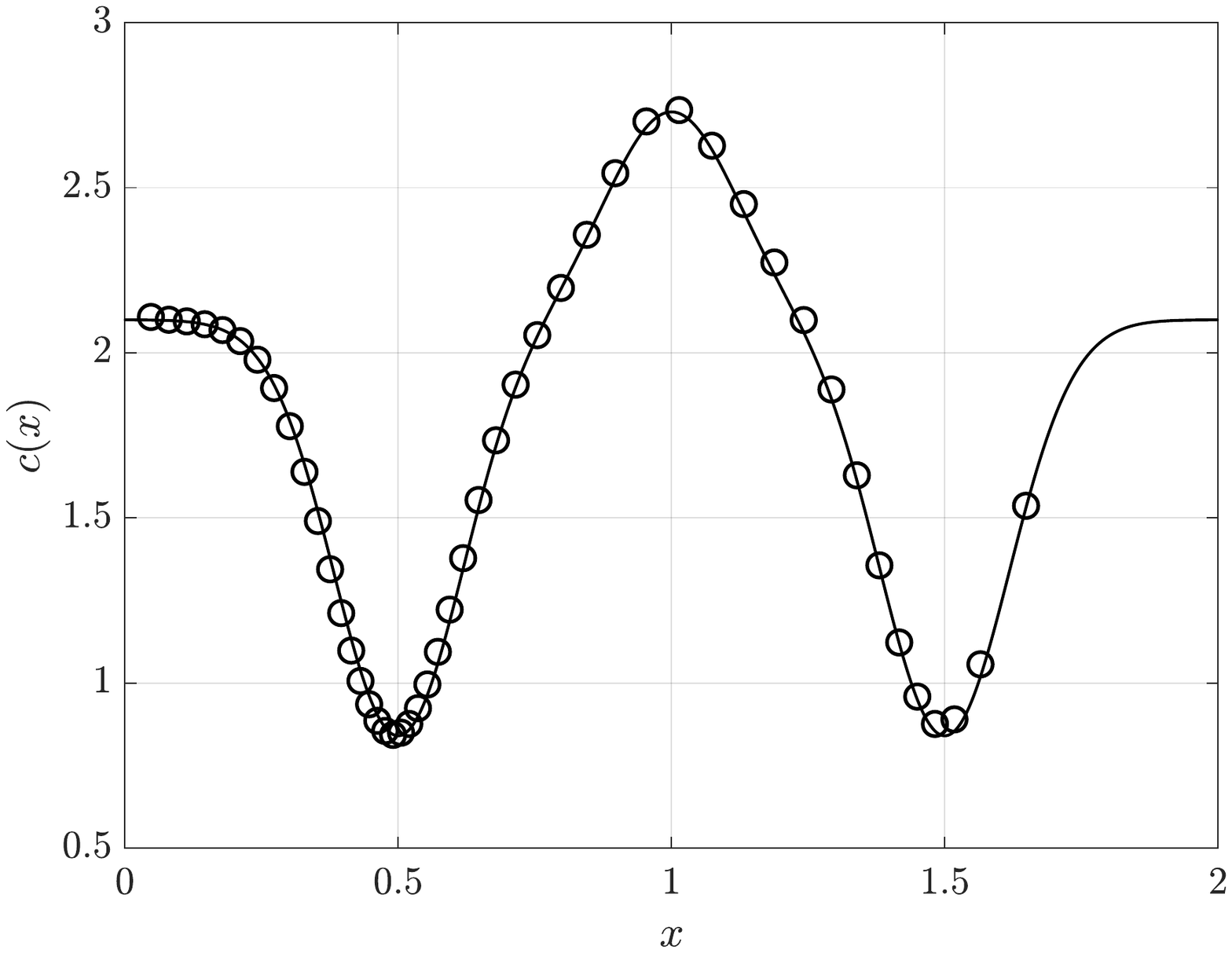}
\caption{Krein embedding reconstruction (circles) of the velocity profile (solid line) based on $n=50$ pole-residue pairs of the impedance function at $x=0$.}
\label{fig:vel_profile1_K_mp}
\end{figure}

\begin{figure}
\centering
\includegraphics[width=0.7\linewidth, trim={0 5cm 0 5cm},clip]{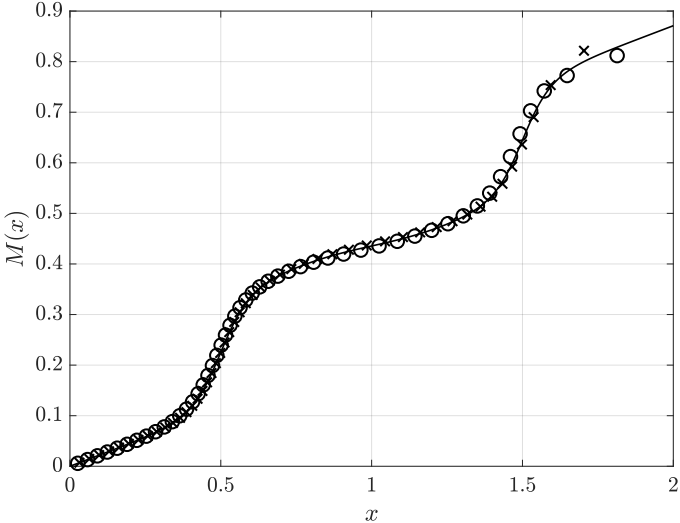}
\caption{Optimal grid embedding reconstruction of the mass function (solid line) based on $n=50$ pole-residue pairs of the impedance function at $x=0$. Crosses: reconstructed mass function values at primary nodes, circles: reconstructed mass function values at dual nodes.}
\label{fig:Mfunction_og_mp}
\end{figure}

\begin{figure}
\centering
\includegraphics[width=0.7\linewidth, trim={2cm 7cm 2cm 7cm},clip]{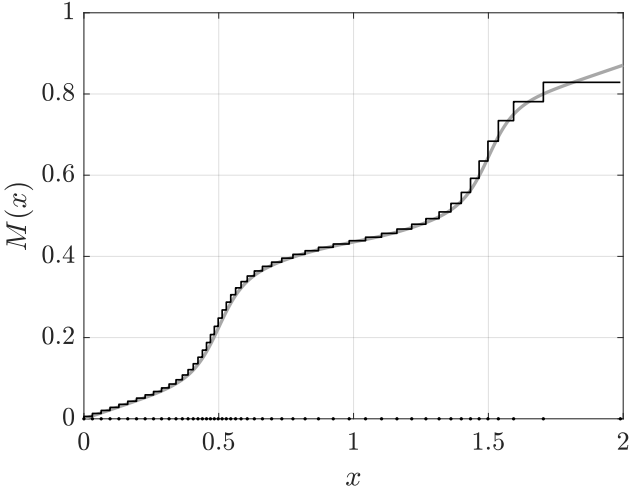}
\caption{Krein embedding reconstruction (black line) of the mass function (grey line) based on $n=50$ pole-residue pairs of the impedance function at $x=0$. The nodes of the Krein grid are shown as solid dots on the $x$-axis.}
\label{fig:Mfunction_K_mp}
\end{figure}

Therefore, let us increase the number of pole-residue pairs to $n=50$. In this case we show the performance of optimal grid and Krein embedding in terms of velocity profile and mass function reconstruction. Figure~\ref{fig:vel_profile1_og_mp} shows the optimal grid embedding reconstruction for $n=50$, and clearly indicates convergence of the optimal grid embedding approach. Figure~\ref{fig:vel_profile1_K_mp} shows the velocity profile reconstruction with the Krein embedding approach, which indicates that Krein embedding converges as well. Finally, mass function reconstructions for optimal grid and Krein embedding are shown in Figures~\ref{fig:Mfunction_og_mp} and \ref{fig:Mfunction_K_mp}, respectively. These reconstruction results indicate that both embedding approaches converge and benefit if more spectral data (pole-residue pairs) of the impedance function is available.    


\begin{figure}
\centering
\includegraphics[width=0.7\linewidth, trim={0 5cm 0 5cm},clip]{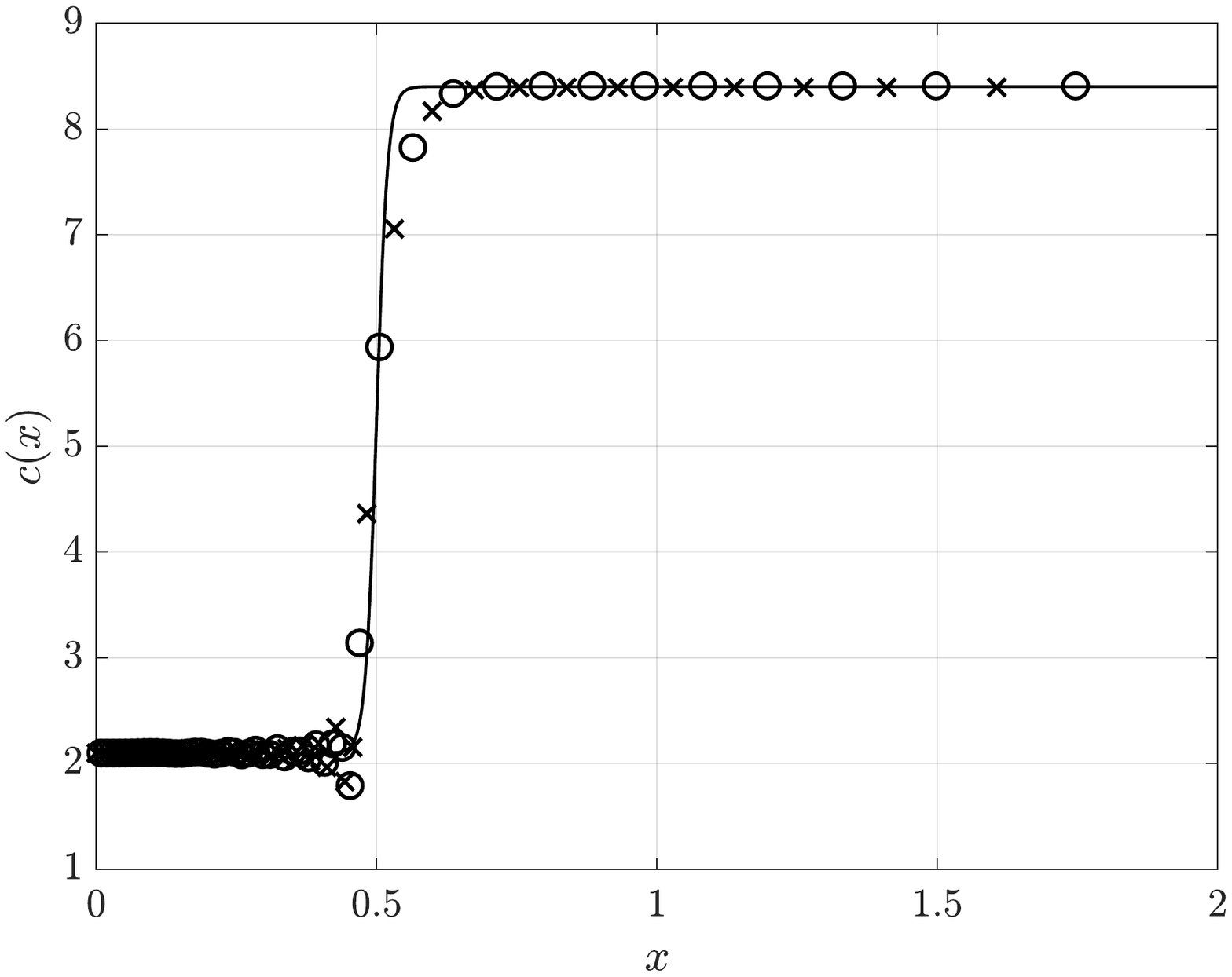}
\caption{Optimal grid embedding reconstruction of a step velocity profile (solid line) based on $n=50$ pole-residue pairs of the impedance function at $x=0$. Crosses: reconstructed velocity values at primary nodes, circles: reconstructed velocity values at dual nodes.}
\label{fig:step_og}
\end{figure}

\begin{figure}
\centering
\includegraphics[width=0.7\linewidth, trim={2cm 7cm 2cm 7cm},clip]{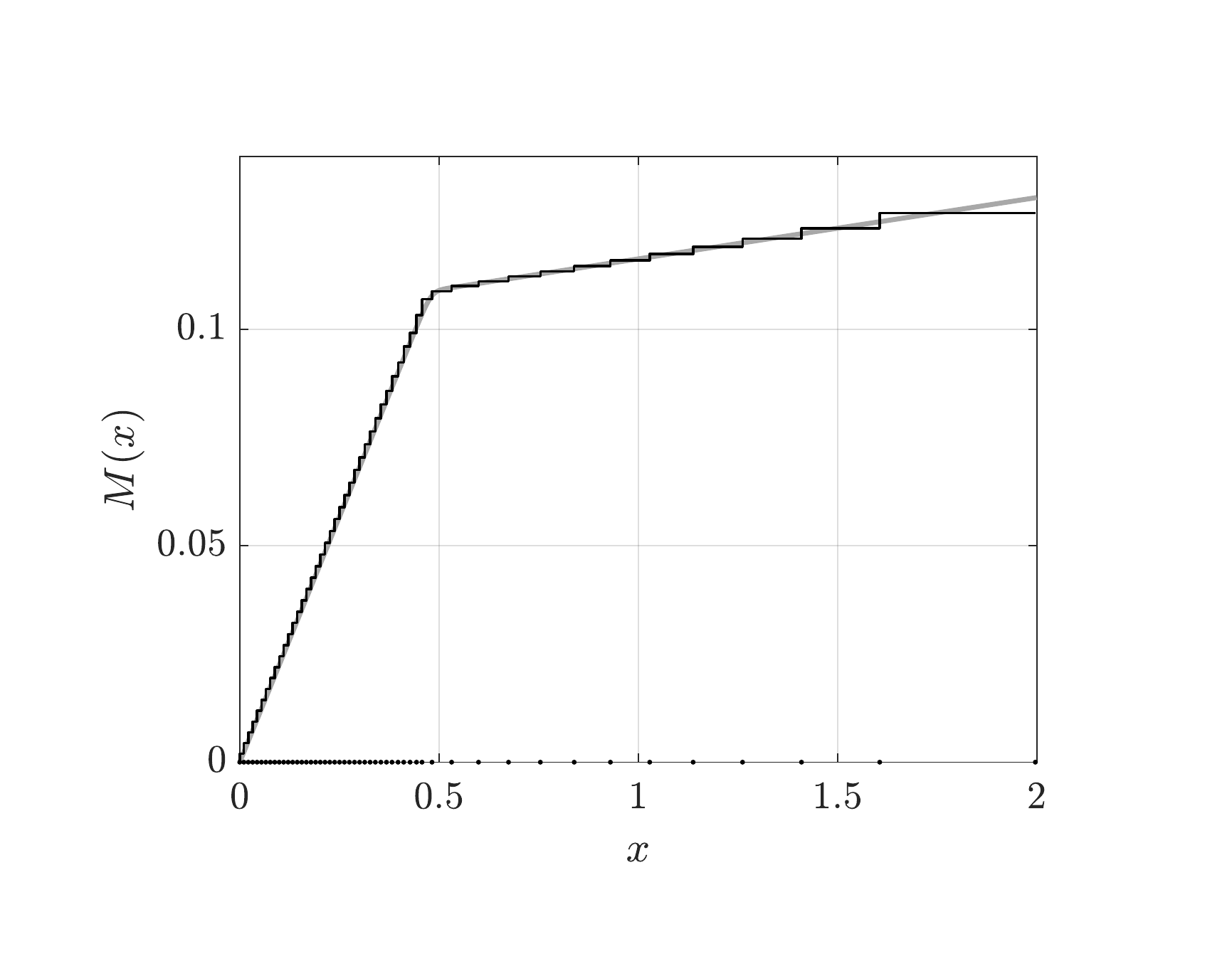}
\caption{Krein embedding reconstruction (black line) of the mass function (grey line) based on $n=50$ pole-residue pairs of the impedance function at $x=0$. The nodes of the Krein grid are shown as solid dots on the $x$-axis. }
\label{fig:step_K}
\end{figure}

Finally, we show how the two embedding approaches reconstruct a smoothed step (piecewise constant) velocity profile. The dashed line in Figure~\ref{fig:step_og} shows the velocity profile along with the optimal grid embedding reconstruction for $n=50$. Clearly, optimal grid embedding captures the exact wave speed away from the step and shows some Gibbs ringing around the step. The Krein embedding reconstruction for $n=50$ of the corresponding mass function is shown in Figure~\ref{fig:step_K}. Here Gibbs ringing is not observed, since an integrated velocity profile is reconstructed. 


\begin{figure}
\centering
\includegraphics[width=0.9\linewidth]{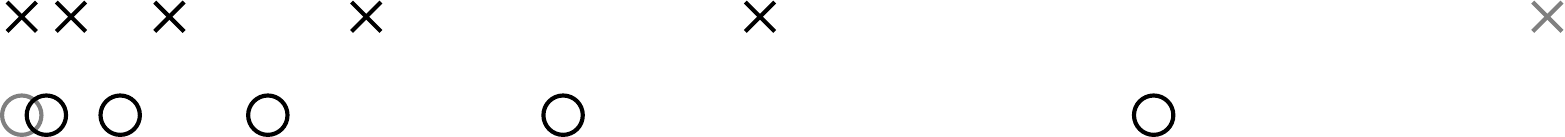}
\caption{Illustration of primary (crosses) and dual (circles) grid nodes of an optimal/Krein grid for a semi-infinite-domain. The grid steps become progressively large as we move away from the measurement point, which coincides with the left-most primary grid node.} 
\label{fig:large_steps}
\end{figure}

\section{Embedding of ROMs on semi-infinite domains}
\label{sec:absorbing}
Up till now, we have considered scattering problems on bounded domains. For semi-infinite domains, that is, domains that extend to infinity in one direction, the optimal and Krein grids may be applied, but these approaches produce excessively large step sizes \cite{Ingerman_etal} as illustrated in Figure~\ref{fig:large_steps}. 

This problem can be avoided, however, by constructing ROMs that can be interpreted as two-point finite-difference discretisations of dissipative first-order continuous wave equations. We refer to this approach as \emph{Krein-Nudelman embedding}~\cite{Krein&Nudelman}.      

To be specific, consider the wave equation in the Laplace domain on a semi-infinite domain, where the wave speed variations are supported on $(0,L)$ and the wave speed is constant $c(x)=c(L)$ for all $x>L$
\begin{equation}
\label{eq:openWaveeq} \hspace{-2cm}
\frac{\mathrm{d}^2u(x,s)}{\mathrm{d}x^2}-  s^2 \frac{1}{c^2(x)}u(x,s)=0, \quad \frac{\mathrm{d}u}{\mathrm{d}x}\Big|_{x=0}=-s,\quad \mbox{and} \quad \lim_{x\rightarrow L}\frac{\mathrm{d}u}{\mathrm{d}x}=\frac{-s}{c(L)}u(L,s).
\end{equation}
The boundary condition at $x=L$ is the Sommerfeld radiation condition only satisfied by waves traveling along the positive $x$-direction (${\rm exp}[-sx/c(L)]$). The inverse problem is the recovery of $M(x)=\int_0^x c(\xi)^{-2}\mathrm{d} \xi$ from impedance measurements $f(s)=u(0,s)$.

The impedance function measured at $x=0$ for an open domain does not have a point spectrum. To implement Krein-Nudelman embedding, we fit the impedance function with rational functions with $n$ complex-conjugate pole-residue pairs. Since the impedance function is passive, the poles and residues have a positive real part. In the presented results, we use vector fit \cite{Gustavsen&Semlyen} to obtain the poles and residues from the impedance function.

With complex poles and residues the diagonal entries (called $\alpha$ in Algorithm~\ref{alg:algo1}) of the tridiagonal matrix obtained from the Lanczos algorithm are no longer zero. Since similar matrices have the same trace we can equate the sum of the eigenvalues in $\bLambda$ to the sum of the diagonal elements of $\bT^{\rm TS}$, called $\alpha_i$'s in the Lanczos algorithm
\begin{equation}
0<{\rm trace}{(\bLambda)} = {\rm trace}{(\bT^{\rm TS})}=\sum_{i=1}^n \lambda_i + \overline{\lambda_i} = \sum_{i=1}^{2n}\alpha_i
\end{equation}
and see that the diagonal entries cannot all vanish. These diagonal elements are interpreted as a loss term in Krein-Nudelman embedding and absorb energy. To facilitate an embedding that allows for absorption, we introduce the complex symmetry tridiagonal matrix  $\bT^{\rm TS;open}$
\begin{equation}
\hspace{-2cm}
\label{eq:TgenericLOSS}
\bT^{\rm TS;open}=
\left(
\begin{array}{ccccc}
\alpha_1 &-\rm{i} (\gamma_1 \hat\gamma_1)^{-\half} & 0 & &0 \\
-\rm{i}(\gamma_1 \hat\gamma_1)^{-\half}  & \alpha_2 & \rm{i}(\gamma_{1} \hat\gamma_{2})^{-\half} &0& \\
\vdots &  \rm{i}(\gamma_{1} \hat\gamma_{2})^{-\half} & \ddots & \ddots\\
\vdots&     & \ddots & \alpha_{2n-1} & -\rm{i}(\gamma_{n} \hat\gamma_{n})^{-\half} \\
0 & \cdots &  0 & -\rm{i}(\gamma_{n} \hat\gamma_{n})^{-\half} &  \alpha_{2n} \\
\end{array}
\right),
\end{equation}
which we obtain from running the Lanczos algorithm with the complex poles and residues obtained from fitting the impedance function of a semi-infinite domain. 

\begin{figure}
\centering
\includegraphics[width=0.9\linewidth]{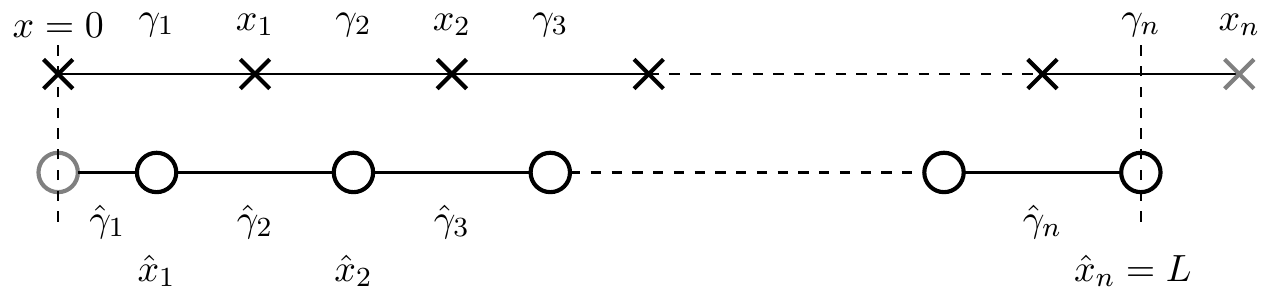}
\caption{Grid used to interpret the ROM as a finite-difference discretisation of the underlying differential operator on a  semi-infinite domain. The crosses represent primary grid nodes and the circles dual grid nodes. Note that the open boundary condition is applied at the last dual node $x=L$.}
\label{fig:grid_open}
\end{figure}

To interpret this tridiagonal matrix using Krein-Nudelman embedding consider a finite-difference discretization of equation~(\ref{eq:openWaveeq}) on the grid displayed in Figure~\ref{fig:grid_open}
\begin{eqnarray}
\label{eq:GenericFD1_3_open}
\frac{\hat u_{j+1}- \hat u_{j}}{c_{j+1}^{-2} \hat\delta_{j+1}} + s {u}_{j} &= 0 \quad \forall j=0,\dots,n-1 \nonumber\\
\frac{u_{j}-  u_{j-1}}{\delta_{j}} + s \hat{u}_{j} &= 0 \quad \forall j=1,\dots,n-1,\\
{-u_{n}}+ c(L)\hat{u}_n+ s\frac{1}{2}\delta_n \hat{u}_{j} &= 0 \nonumber
\end{eqnarray}
and compare it to finite difference discretization introduced by the tridiagonal matrix in equation~(\ref{eq:TgenericLOSS})
\begin{eqnarray}
\label{eq:GenericFD1_3_3}
\frac{\hat u_{j+1}-\hat u_{j}}{\hat\gamma_{j+1}} + (\alpha_{2j+1}+s) {u}_{j} &= 0, \qquad \forall j=0,\dots,n-1, \nonumber \\
\frac{u_{j}-u_{j-1}}{\gamma_{j}} + (\alpha_{2j}+s) \hat{u}_{j} &= 0, \qquad \forall j=1,\dots,n.
\end{eqnarray}
Note that the Sommerfeld condition is applied on the last dual grid point, whereas the Dirichlet condition in the bounded domain was applied on the last primary grid node.

The diagonal of $\bT^{\rm TS;open}$ appears as a loss terms in the equation, and can thus be used to define integrated primary and dual loss coefficients
\begin{equation}
r_i=\hat\gamma_i \alpha_{2i-1} \quad {\rm and } \quad \hat{r}_i=\gamma_i \alpha_{2i}.
\end{equation}
A Sommerfeld radiation condition would be equivalent to all $r_i$ and $\hat{r}_i$ values vanishing except for $\hat{r}_n$.

\begin{figure}
\centering
\includegraphics[width=0.7\linewidth, trim={2cm 7cm 2cm 7cm},clip]{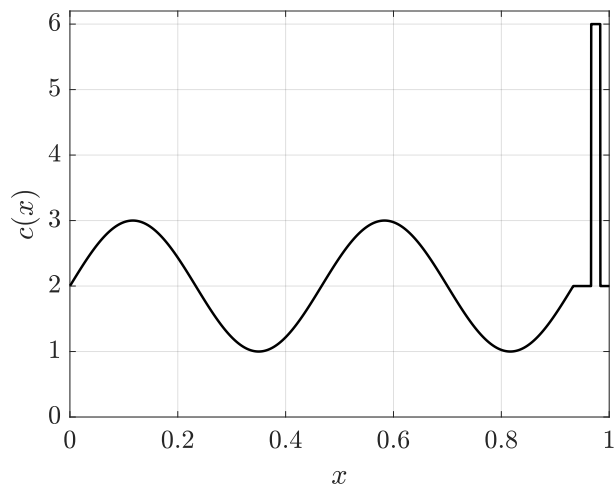}
\caption{Wave speed model for the semi-infinite domain example. A reflector is placed at the end of the medium.}
\label{fig:OpenSin6_Medium}
\end{figure}

To illustrate this embedding we consider the inverse problem of recovering the velocity profile in Figure~\ref{fig:OpenSin6_Medium} from impedance measurements at $x=0$. We use $n=121$ pole-residue pairs to fit the impedance function and compute $\bT^{\rm TS;open}$ using the complex symmetric Lanczos algorithm. 

\begin{figure}
\centering
\includegraphics[width=0.7\linewidth, trim={2cm 7cm 2cm 7cm},clip]{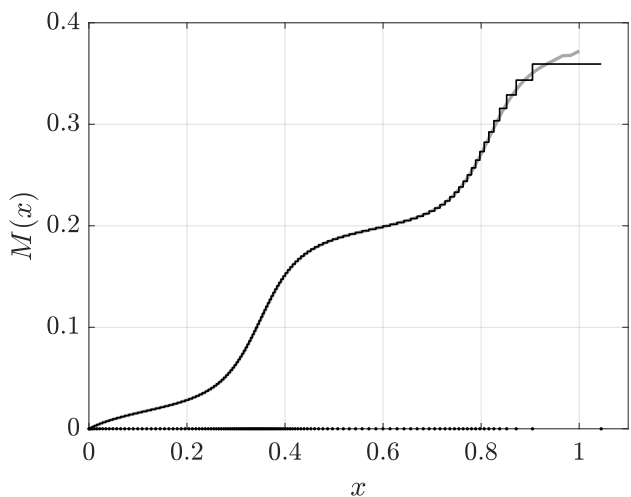}
\caption{Mass function $M(x)$ (solid black line) and Krein embedding grid reconstruction of this function (grey solid line) based on $n=121$ pole-residue pairs approximating the impedance function at $x=0$. }
\label{fig:OpenSin6_Embedding}
\end{figure}

The Krein-Nudelman embedding of the obtained mass function is shown alongside the true mass function in Figure~\ref{fig:OpenSin6_Embedding}. The smooth velocity variations are recovered well, whereas the reflector is not recovered.

\begin{figure}
\centering
\includegraphics[width=0.7\linewidth, trim={2cm 7cm 2cm 7cm},clip]{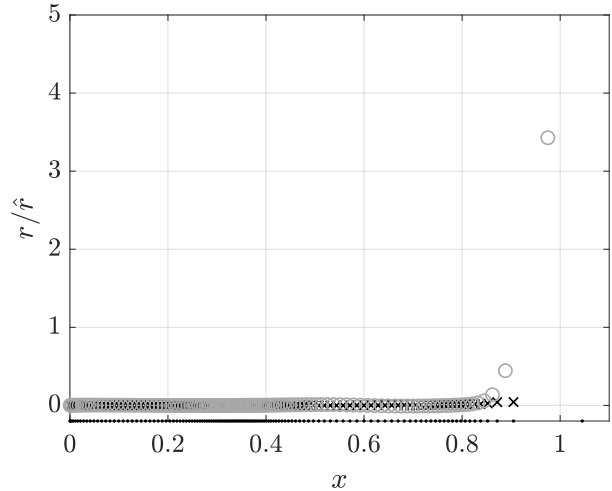}
\caption{Embedding: Primary $r_i$ ($\times$) losses embedded on the Krein grid  $\chi^{\rm Krein}$ and dual losses $\hat{r}_i$ ($\circ$) embedded centered between Krein grid nodes. An effective absorbing inclusion is placed at the end of the domain by the Krein embedding.}
\label{fig:OpenSin6_Loss}
\end{figure}

In Figure~\ref{fig:OpenSin6_Loss} the loss coefficients are embedded into space for visualisation purposes. The primary losses are displayed on the Krein grid and the dual losses in between Krein grid nodes.  Almost all losses $r/\hat{r}$ vanish up until the last reflector in the medium, where the last thee $\hat{r}_i$ contribute to an effective loss term. The recovered $r$ and $\hat{r}$ terms do not exactly correspond to a Sommerfeld radiation condition, but rather an effective absorbing medium being placed at the end of the Krein embedding. Nonetheless, the velocity profile up to the last reflector is well recovered at which point the loss coefficients increase. Note that the last grid point of the Krein embedding is outside the interval of wave speed variations, since the Sommerfeld condition is prescribed on a dual grid node in between the Krein-grid nodes.

Here we reach a fundamental limitation of the Krein embedding. From the fitted pole residue form, the Krein embedding cannot distinguish between an infinite lossless medium, which does not have a point spectrum, and a finite-absorbing medium. See the Appendix for a further discussion. 

However, up until the Krein embedding recovers an  ``effective'' absorbing medium the embedding can be used to solve inverse problems. This is in some sense physical as one typically does not know from which point the infinite domain has constant medium parameters. Compared to an optimal grid approach, the Krein-Nudelman approach only produces a grid in the interval where the variations of the wave speed are supported, whereas the optimal grid approach places grid points beyond the variations of the medium, which is not useful for inversion.

\section{Conclusions} 
\label{sec:concl}
In this paper, we discussed two so-called ROM embedding procedures to solve inverse scattering problems based on spectral impedance data collected at one end of an interval of interest. More precisely, we have shown that starting from impedance data in pole-residue form, the Lanczos algorithm can be used to construct ROMs that can be interpreted as two-point finite-difference discretisations of first-order wave equations. We call this ROM embedding and presented two such embedding procedures, namely, optimal grid and Krein embedding. 

The optimal grid embedding procedure is based on the theory of optimal grids and was implemented in terms of spatial coordinates instead of travel time coordinates. In this approach, the wave speed profile of an inhomogeneous medium can be reconstructed if we first determine an optimal reference grid for a (homogeneous) reference medium. The use of a reference grid is avoided, however, in the Krein embedding approach, and convergence results obtained for optimal grid embedding can be used to demonstrate the convergence of Krein embedding as well. Furthermore, since no reference grid is necessary for Krein embedding, it may be possible to apply this procedure to semi-infinite domains. Here we run into uniqueness problems, however, since with an impedance function in pole-residue form, it is impossible to distinguish between a lossy medium on a bounded interval and a lossless medium on a semi-infinite interval. Nevertheless, numerical examples indicate that the wave speed profile is reconstructed on a semi-infinite domain up till the last reflector that is present within this domain. Future work will focus on the development of a ROM embedding procedure for semi-infinite domain problems.

\section*{Appendix}
Consider the system $Au=\lambda u +g$, where $A$ is a Hermitian nonnegative operator, and $u$ and $g$ are elements belonging to a Hilbert space. We introduce the transfer function of this system as $f(\lambda) = g^\ast u$, where the asterisk denotes the Hermitian inner product on the Hilbert space. This transfer function is a Stieltjes function. Recall that a function $f(\lambda)$ is a Stieltjes function if  
\begin{equation}
\label{eq:S}
f(\lambda) = 
\int_{z=-\infty}^0 
\frac{\rho(z)}{\lambda -z} \, \mbox{d}z, \quad \lambda \in \mathbb{C}\setminus\mathbb{R}_-,
\end{equation}
for some positive distribution $\rho(z)$ defined as the generalised derivative of a probability measure.

Subsequently, we consider the first-order system 
\[
\left(
\begin{array}{cc}
-\sigma & L \\
L^\ast & 0
\end{array}
\right)
\left(
\begin{array}{c}
\tilde{u} \\
v
\end{array}
\right)
= 
s 
\left(
\begin{array}{c}
\tilde{u} \\
v
\end{array}
\right) +
\left(
\begin{array}{c}
g \\
0
\end{array}
\right),
\]
on a Hilbert space, where $L^\ast$ is the adjoint of operator $L$ and $\sigma \geq 0$. The transfer function of this system is introduced as $h(s) = g^\ast \tilde{u}$. A lossless first-order system is characterised by $\sigma=0$ everywhere, while if $\sigma>0$ in some domain in space, we refer to the above system as a lossy first-order system. 

Eliminating $v$ in the lossless case, we obtain the Hermitian system $A\tilde{u}=\lambda \tilde{u} + sg$ with $A=LL^\ast$ and $\lambda=s^2$. Here, the transfer functions of the Hermitian and first-order system are related by $f(\lambda) = \lambda^{-1/2} h(\lambda^{1/2})$ and $h$ satisfies the three criteria for a passive SISO system, that is, we have 
\begin{enumerate}
\item $h(s)$ has no poles in $\mathbb{C}_+$,
\item $h(\bar{s}) = \bar{h}(s)$ for every $s \in \mathbb{C}_+$, where the overbar denotes complex conjugation, and 
\item $\mbox{Re}[h(\mbox{i}\omega +0)] \geq 0$ for $\omega \in \mathbb{R}$.
\end{enumerate}

For a lossy first-order system, the transfer function $h(s)$ still satisfies the above three criteria, but we no longer arrive at a Hermitian form when $v$ is eliminated from the first-order system. However, we can establish an isomorphism between passive transfer functions and Stieltjes functions. More precisely, we have 

\bigskip
\noindent
{\sc Proposition~2}. Using $f(\lambda) = \lambda^{-1/2} h(\lambda^{1/2})$, every Stieltjes function $f(\lambda)$ can be transformed to a passive transfer function $h(\lambda^{1/2})$ and vice versa.

\bigskip
\noindent
{\it Proof:} Consider the transformation $f \rightarrow h$. First note that from the definition (\ref{eq:S}) of a Stieltjes function it follows that $f(\lambda)$ is an analytic function of $\lambda$ on $\mathbb{C}\setminus \mathbb{R}_-$. Setting $s=\lambda^{1/2}$, it follows from $h(s)=sf(s^2)$ and the analyticity of $f$ that $h(s)$ is analytic in $\mathbb{C}_+$ and therefore $h(s)$ has no poles in $\mathbb{C}_+$. Furthermore, since $\rho$ is real (positive), we have $\bar{f}(s^2) = f(\bar{s}^2)$ and it follows that $\bar{h}(s) = h(\bar{s})$ with $s \in \mathbb{C}_+$. Finally, we have 
\begin{equation}
\lim_{\epsilon \downarrow 0}
\mbox{Re}[h(\epsilon+ \mbox{i}\omega)] = 
\lim_{\epsilon \downarrow 0}
\frac{1}{2} 
\left[
h(\epsilon+ \mbox{i}\omega) + h(\epsilon-\mbox{i}\omega) 
\right] = \pi |\omega| \rho(-\omega^2)>0,
\end{equation}
where we have used $h(s)=sf(s^2)$ and 
\[
\rho(\lambda) = 
\frac{f(\lambda +\mbox{i}0) - f(\lambda - \mbox{i}0)}{2\pi\mbox{i}}, \quad 
\mbox{with $\lambda \in \mathbb{R}_-$.}
\]
Reasoning in reverse order gives the ``vice versa" result. \hfill $\Box$

\bigskip
\noindent
Now to each Stieltjes function there corresponds a unique Stieltjes-Krein string \cite{Kac&Krein} and a consequence of the above proposition is that we cannot distinguish between lossy media and lossless media on a semi-infinite domain, since in both cases the transfer function is passive.         
%
%

\section*{Acknowledgements}
We thank Liliana Borcea, Alex Mamonov, and Mikhail Zaslavsky for many stimulating discussions. The work of Vladimir Druskin was financially supported by AFOSR grants FA 955020-1-0079 and FA9550-23-1-0220, and NSF grant DMS-2110773. The work of Elena Cherkaev was financially supported by NSF grant DMS-2111117. The work of Murthy Guddati was financially supported by NSF grant DMS-2111234. The work of J\"{o}rn Zimmerling was financially supported by NSF grant DMS-2110265. This support is gratefully acknowledged. 

\end{document}